\documentclass[12]{article} 
\usepackage{amsfonts,latexsym}
\usepackage{theorem}
\usepackage{graphicx}
\usepackage{amsmath}
\usepackage[a4paper, margin=2.5cm]{geometry}

\def\CC{{\mathbb C}}
\def\DD{{\mathbb D}}
\def\NN{{\mathbb N}}

\def\RR{{\mathbb R}}
\def\TT{{\mathbb T}}

\newtheorem{theorem}{Theorem}
\newtheorem{lemma}{Lemma}

\newcommand{\boite}{\mbox{} \hfill \mbox{\rule{2mm}{2mm}}}

\title{Constrained $L^2$-approximation by polynomials\\ 
on subsets of the circle}
\author{L. Baratchart\footnote{INRIA, BP 93, 
06902 Sophia-Antipolis Cedex, FRANCE}, \and
J. Leblond\footnotemark[1], \and F. Seyfert\footnotemark[1]}  

\begin{document}

\maketitle
\hfill
\emph{To the memory of Andr\'e Boivin}
\section{Abstract}
We study best approximation to a given function, in the least square sense 
on a subset of the unit circle, by polynomials of given degree
which are pointwise bounded on the complementary subset.
We show that the solution to this problem,  as the degree goes large, 
converges to the solution of a bounded  extremal problem for analytic 
functions which is instrumental in system identification. We provide a 
numerical example on real data from a hyperfrequency filter.

\section{Introduction}
\label{sec1}
This paper deals with best approximation to a square summable function,
on a finite union $I$ of arcs of the unit circle $\TT$, by a polynomial of 
fixed degree which is bounded by 1 in modulus on the complementary system
of arcs $J=\TT\setminus I$. This we call, for short, 
the polynomial problem. We are also concerned 
with the natural limiting version when the degree 
goes large, namely
best  
approximation in $L^2(I)$ by a Hardy function of class $H^2$
which is bounded by 1 on $J$. To distinguish this issue from
the polynomial problem, we term it the analytic problem. The latter
is a variant, involving mixed norms, 
of constrained extremal problems for analytic functions considered in
\cite{KN,partII,BLPprep,LP99,schneck1}. As we shall see, solutions to
the polynomial problem
converge to those of the analytic problem as the  degree tends to infinity, 
in a sense to be made precise  below. This is why solving 
for high degree  the polynomial 
problem (which is finite-dimensional)  is an interesting 
way to regularize and approximately solve the analytic problem 
(which is infinite-dimensional). This is the gist of the present work.

Constrained extremal problems for analytic functions, in particular the 
analytic problem defined above,  can be set up more
generally in the 
context of weighted approximation, {\it i.e.} seeking best approximation
 in $L^2(I,w)$ where $w$ is a weight on $I$.
In fact, that kind of generalization is useful for applications 
as we shall see.
As soon as $w$ is invertible in $L^\infty(I)$, though, 
such a weighted problem turns out to be equivalent to
another one with unit weight, hence the present formulation 
warrants  most practical situations. This property allows one to carry  
the analytic problem over to more general curves than the circle.
In particular, in view of the 
isomorphism between Hardy spaces of the disk and the half-plane 
arising by composition with a M\"obius transform \cite[ch. 10]{Duren},
best approximation in $L^2(I)$ from $H^2$ of the disk
can be converted to weighted best approximation 
in $L^2(\mathfrak{I},w)$ from the Hardy space $\mathfrak{h}^2$ of a half-plane
with $\mathfrak{I}$ a finite union of bounded intervals on the line and
$w$ a weight arising from the derivative of the M\"obius transform.
Since this weight is boundedly invertible on $\mathfrak{I}$, it
follows that the analytic problem
on the circle and its analog on the line are equivalent.
One may  also define another Hardy space $\mathcal{H}^2$, say  of the right 
half-plane as the space of analytic functions whose $L^2$-means over
vertical lines are uniformly bounded. Then, best approximation in $L^2(I)$ 
from $H^2$ is equivalent to best approximation from $\mathcal{H}^2$
in $L^2(\mathfrak{I})$, {\it i.e.} weight is no longer needed.
Of course, such considerations hold for many other domains and boundary curves 
than the half-plane and the line, but the 
latter are of special significance to us as we now explain.

Indeed, on the line,
constrained extremal problems for analytic functions naturally 
arise in Engineering when studying
deconvolution  issues, in particular
those pertaining to system identification and design.
This motivation is stressed in 
\cite{KN, partI, aakblp, thesefab, BLPprep}, whose
results are effectively used  today to identify microwave devices
\cite{RRbglosw,OSM}.  More precisely,
recall that a linear time-invariant 
dynamical system is just a convolution operator, 
hence the Fourier-Laplace transform
of its output is that of its input times 
the Fourier-Laplace transform of its kernel. The latter is called the 
transfer-function. Now, by feeding periodic inputs to a stable system, one can 
essentially recover the transfer function pointwise on the line, 
but typically in
a restricted  range of frequencies only, corresponding to 
the passband of the system, say $\mathfrak{I}$ \cite{DFT}. 
Here, the type of stability under consideration  impinges 
on the smoothness of the transfer function as well as on the precise 
kind of recovery
that can be achieved, and we refer the reader to \cite[Appendix 2]{BCS} for a 
more thorough analysis. For the present discussion, it suffices to 
assume that the system is stable in the $L^2$ sense, {\it i.e.} that it
maps square summable inputs to square summable outputs.
Then, its transfer function lies in $H^\infty$ of the half-plane 
\cite{Partingtonbo}, and to identify it
we are led to approximate the measurements   on  
$\mathfrak{I}$ 
by a Hardy function with a bound on its modulus.
Still, on $\mathfrak{I}$, a natural criterion from the stochastic 
viewpoint is $L^2(\mathfrak{I},w)$, where the weight $w$ is the reciprocal of 
the pointwise covariance of the noise assumed to be additive \cite{PRS}.
Since this covariance is boundedly invertible on $I$, we face  
an  analytic problem on the line upon normalizing
the bound on the transfer function to be 1. 
 This stresses how the 
analytic problem on the line, which can be mapped back to the circle, 
connects to system identification. Now, this analytic problem is 
convex but infinite-dimensional.  Moreover, as 
Hardy functions have no discontinuity of the first kind
on the boundary \cite[ch. II, ex. 7]{Garnett} 
and since the solution to an analytic problem generically
has exact modulus 1 on $J$, as we prove later on,  it will typically
oscillate at the endpoints of $I$, $J$ which is unsuited.
One way around these difficulties is to solve the polynomial problem
for sufficiently high degree, 
as a means to 
regularize and approximately solve the analytic one. This was an initial
motivation by the authors to write the present  paper, and we 
provide the reader in Section \ref{numerique} with a numerical example on real data 
from a hyperfrequency filter. It must be said that the polynomial 
problem itself has numerical issues: though  it is convex in
finitely many variables,
bounding the modulus on $J$  involves
infinitely many convex constraints
which  makes it of so-called semi-infinite programming type.
A popular technique to handle such problems is through linear matrix 
inequalities,  but we found it easier to approximate from below
the polynomial 
problem by a finite-dimensional one with finitely many constraints,
in a demonstrably convergent manner as the number of these constraints gets large.

The organization of the paper is as follows. In section \ref{notprel} we set
some notation and we recall standard properties of Hardy spaces.
We state the polynomial and analytic problems in Section \ref{wposed},
where we also show  they are well-posed. Section \ref{critpsec} deals with
the critical point equations characterizing the solutions, and
with  convergence of the polynomial problem to the analytic one.
Finally, we report on some numerical experiment in Section \ref{numerique}.

\section{Notations and preliminaries}
\label{notprel}
Throughout we let $\TT$ be the unit circle 
and $I\subset \TT$ a finite union of nonempty open  arcs
whose complement $J=\TT\setminus I$ has nonempty interior. 
If $h_1$ (resp. $h_2$) is a function defined on a set containing $I$ 
(resp. $J$), we put $h_1\vee h_2$ 
for the concatenated function, defined on the whole of $\TT$, which is $h_1$
on $I$ and $h_2$ on $J$. 

For $E\subset\TT$, we let $\partial E$ and 
$\stackrel{\circ}{E}$ denote respectively the boundary and the interior
of $E$ when viewed as a subset of $\TT$; we also let $\chi_E$ for the
characteristic function of $E$ and
$h_{|_E}$ for the restriction of $h$ to $E$. Lebesgue measure on $\TT$ is just the image of Lebesgue measure on
$[0,2\pi)$ under the parametrization $\theta\mapsto e^{i\theta}$.
We denote by $|E|$ the measure of a measurable subset $E\subset \TT$,
and if $1\leq p\leq\infty$ we write $L^p(E)$ 
for the familiar Lebesgue space of
(equivalence classes of a.e. coinciding) complex-valued measurable functions 
on $E$ with norm
\[\|f\|_{L^p(E)}=\left(\frac{1}{2\pi}
\int_E|f(e^{i\theta})|^p\,d\theta\right)^{1/p}<\infty
\quad \text{if }1\leq p<\infty, \qquad
\|f\|_{L^\infty(E)}={\rm ess.}\sup_{\theta \in 
E}|f(e^{i\theta})|<\infty.\]
We sometimes indicate by $L^p_{\RR}(E)$ the real subspace of
real-valued functions. 
We also set
\begin{equation}
\label{defdual}
\langle f,g \rangle_E= \frac{1}{2\pi}\int_E f(e^{i\theta})\overline{g(e^{i\theta})}\,
d\theta
\end{equation}
whenever $f\in L^p(E)$ and $g\in L^q(E)$ with $1/p+1/q=1$. If $f$ and
$g$ are defined on a set containing $E$, we  write for simplicity
$\langle f,g\rangle_E$ to mean $<f_{|_E},g_{|_E}>$ and
$\|f\|_{L^p(E)}$ to
mean $\|f_{|_E}\|_{L^p(E)}$.
Hereafter 
$C(E)$ stands for the space of 
bounded complex-valued continuous functions on
$E$ endowed with the $\sup$ norm, while $C_\RR(E)$ indicates real-valued 
continuous functions. 

Recall that the Hardy space $H^p$ is the closed
subspace of $L^p(\TT)$ consisting of functions whose Fourier coefficients
of strictly negative index do vanish. 
We refer the reader to \cite{Garnett} for standard facts on Hardy spaces,
in particular those recorded  hereafter.
Hardy functions are the nontangential limits a.e.
on $\TT$ of functions  holomorphic in the unit disk $\DD$ 
having uniformly bounded $L^p$ 
means over all circles centered at $0$ of radius less than 1:
\begin{equation}
\label{defintHp}\|f\|_{H^p}=\sup_{0\leq r<1}\left(\frac{1}{2\pi}\int_0^{2\pi} 
|f(re^{i\theta})|^p\,d\theta\right)^{1/p}
\quad \text{if }1\leq p<\infty, \qquad
\|f\|_{H^\infty}=\sup_{z\in\DD}|f(z)|.
\end{equation}
The correspondence between such a holomorphic function $f$ and its non 
tangential limit $f^\sharp$ is one-to-one and even isometric,
namely the supremum in \eqref{defintHp} is equal to
$\|f^\sharp\|_p$, thereby allowing us to identify $f$ and $f^\sharp$ and 
to drop the superscript $\sharp$.
Under this 
identification, we  regard members of $H^p$ both as functions in
$L^p(\TT)$ and as holomorphic
functions in the variable $z\in\DD$, but the argument (which belongs 
to $\TT$ 
in the former  case and to $\DD$ in the latter) helps preventing confusion.
It holds in fact that $f_r(e^{i\theta})=f(re^{i\theta})$ converges as
$r\to 1^-$ to
$f(e^{i\theta})$ in $L^p(\TT)$ when $f\in H^p$ and $1\leq p<\infty$.
It follows immediately from \eqref{defintHp} and H\"older's inequality that,
whenever 
$g_1\in H^{p_1}$ and $g_2\in H^{p_2}$, we have $g_1g_2\in H^{p_3}$ 
if $1/p_1+1/p_2=1/p_3$.

Given $f\in H^p$, its values on $\DD$ are obtained from its values on $\TT$
through a Cauchy as well as a Poisson integral \cite[ch. 17, thm 11]{Rudin}, 
namely:
\begin{equation}
\label{Cauchy1}
f(z) = \frac{1}{2\,i \, \pi} \, \int_{\TT}
\frac{f(\xi)}{\xi-z}\,d\xi\,,~~~~{\rm and~also}~~~~
f(z)=\frac{1}{2 \, \pi} \, \int_{\TT} \mbox{\rm Re}\left\{\frac{e^{i\theta} + 
z}{e^{i\theta} - z} \right\}\, f(e^{i\theta})\, d\theta \, ,~~ ~~z\in\DD,
\end{equation}
where the right hand side of
the first equality in \eqref{Cauchy1} is a line integral.
The latter immediately implies  that the Fourier coefficients of a Hardy function on the circle are the Taylor coefficients of its power series expansion at
0 when viewed as a holomorphic function on $\DD$. In this connection,
the space $H^2$ is especially simple to describe:
it consists of those holomorphic functions $g$ in $\DD$ whose 
Taylor coefficients at 0 are square summable, namely
\begin{equation}
 g(z)=\sum_{k=0}^\infty a_kz^k:\quad
\|g\|_{H^2}^2:=\sum_{k=0}^\infty |a_k|^2<+\infty,\qquad
g(e^{i\theta})=\sum_{k=0}^\infty a_ke^{ik\theta},
\end{equation}
where the convergence of the last Fourier series holds in $L^2(\TT)$ by 
Parseval's theorem (and also pointwise a.e. by Carleson's theorem but 
we do not need this deep result).
Incidentally, let us mention that for no other value of
$p$ is it known how to characterize
$H^p$ in terms of the size of its Fourier coefficients.

By  the Poisson representation ({\it i.e.} the second integral in \eqref{Cauchy1}), a Hardy
function $g$ is also uniquely represented, up to a
purely imaginary constant, by its real part $h$ on $\TT$ according to:
\begin{equation}
g(z) = i\mbox{\rm Im}g(0)\,+\,\frac{1}{2 \, \pi} \, \int_{\TT} \frac{e^{i\theta} + 
z}{e^{i\theta} - z} \, h(e^{i\theta})\, d\theta \, ,~~ ~~z\in\DD.
\label{RHt}
\end{equation}
The integral in (\ref{RHt}) is called the
\emph{Riesz-Herglotz transform} of $h$ and, whenever $h\in L_{\RR}^1(\TT)$, it
defines a holomorphic function in $\DD$ which is real at 0 and whose
nontangential limit exists a.e. on $\TT$ with real part equal to $h$.
Hence 
the Riesz-Herglotz transform (\ref{RHt}) assumes the form 
$h(e^{i\theta})+i\widetilde{h}(e^{i\theta})$ a.e. on $\TT$, where the
real-valued function $\widetilde{h}$ is said to be \emph{conjugate} to $h$.
It is a theorem of M. Riesz
\cite[chap. III, thm 2.3]{Garnett}
that  if $1<p<\infty$, then $\widetilde{h}\in L_{\RR}^p(\TT)$
when $h\in L_{\RR}^p(\TT)$. This neither holds  
for $p=1$ nor for $p=\infty$. 

A nonzero $f\in H^p$ can be uniquely
factored as $f=jw$ where
\begin{equation}
\label{def-outer}
w(z)=\exp\left\{\frac{1}{2\pi}\int_0^{2\pi}\frac{e ^{i\theta}+z}{e
    ^{i\theta}-z}\log|f(e ^{i\theta})|\, d\theta\right\}
\end{equation}
belongs to $H^p$ and is called the {\em outer factor} of $f$, while
$j\in H^\infty$ has modulus 1 a.e. on $\TT$ and is called the
{\em inner factor} of $f$. 
That $w(z)$ in (\ref{def-outer}) is well-defined rests on the fact that
$\log|f|\in L^1$ if $f\in H^1\setminus\{0\}$; it
entails that a $H^p$ function cannot vanish 
on a subset of strictly positive Lebesgue measure on $\TT$ unless it is 
identically 
zero. For simplicity, we often say that a function is outer (resp. inner) 
if it is equal, up to a unimodular multiplicative constant, to its outer (resp. inner) factor.

Closely connected to Hardy spaces is the Nevanlinna class $N^+$, consisting
of holomorphic functions in $\DD$
that can be factored as $jE$, where $j$ is an inner function and $E$ an outer 
function of the form
\begin{equation}
\label{outa}
E(z)=\exp\left\{\frac{1}{2\pi}\int_0^{2\pi}\frac{e ^{i\theta}+z}{e
    ^{i\theta}-z}\log \rho(e ^{i\theta})\, d\theta\right\}\,,
\end{equation}
with $\rho$ a positive function
such that $\log\rho\in L^1(\TT)$ (though $\rho$ itself may not
be summable). Such a function 
has nontangential limits of modulus $\rho$ a.e. on $\TT$.
The Nevanlinna class is instrumental  in that
$N^+\cap L^p(\TT)=H^p$, see  \cite[thm 2.11]{Duren} or 
\cite[5.8, ch.II]{Garnett}. Thus, formula (\ref{outa}) defines a $H^p$-function
if and only if $\rho\in L^p(\TT)$.

Let $\hat{\CC}=\CC \cup \{\infty\}$ be the Riemann sphere. The Hardy space ${\bar H}^p$ of $\hat{\CC}\setminus\overline{\DD}$ 
can be given a treatment parallel to $H^p$ upon changing $z$ into $1/z$.
Specifically,
${\bar H}^p$  consists of  functions in $L^p(\TT)$ 
whose Fourier coefficients of strictly positive index do vanish; these
are, a.e. on $\TT$, the complex conjugates of $H^p$-functions, and
they can also be viewed as nontangential limits of 
functions analytic in $\hat{\CC}\setminus \overline{\DD}$
having uniformly bounded $L^p$ means 
over all circles centered at $0$ of radius bigger than $1$.
We further single out the subspace ${\bar H}_0^p$ of ${\bar H}^p$,
consisting of functions vanishing at infinity or, equivalently, 
having vanishing mean on $\TT$. Thus, a function belongs to  ${\bar H}_0^p$
if, and only if it is of the form 
$e^{-i\theta}\overline{g(e^{i\theta})}$ for some $g\in H^p$. For
$G\in{\bar H}^p_0$, the Cauchy formula assumes the form~:
\begin{equation}
\label{Cauchy2}
G(z) = 
\frac{1}{2\,i \, \pi} \, \int_{\TT}
\frac{G(\xi)}{z-\xi}\,d\xi\,,~~~~z\in\hat{\CC}\setminus\overline{\DD}.
\end{equation}

It follows at once from the Cauchy formula that
the duality product $\langle ~,~\rangle_{\TT}$ makes $H^p$ and ${\bar H}_0^q$ 
orthogonal to each other, and it reduces to the familiar scalar product 
when $p=q=2$. In particular, we have 
the orthogonal decomposition~:
\begin{equation}
\label{orthogdec}
L^2(\TT)=H^2\oplus{\bar H}_0^2.
\end{equation}
For $f\in C(\TT)$ and $\nu\in\mathcal{M}$, the space of complex 
Borel measures on 
$\TT$, we set 
\begin{equation}
\label{pairdual}
\nu{\bf .} f=\int_\TT f(e^{i\theta})\,d\nu(\theta)
\end{equation}
and this pairing induces an isometric 
isomorphism between $\mathcal{M}$ (endowed with the norm of the
total variation) and 
the dual of $C(\TT)$ \cite[thm 6.19]{Rudin}. If we let
$\mathcal{A}\subset H^\infty$ designate the disk algebra of
functions analytic in $\DD$ and continuous on $\overline{\DD}$, 
and if $\mathcal{A}_0$ indicates those functions in $\mathcal{A}$
vanishing at zero, it is easy to
see that $\mathcal{A}_0$ is the orthogonal space under
(\ref{pairdual}) to those measures whose Fourier
coefficients of strictly negative index do vanish. 
Now, it is a fundamental theorem 
of F. and M. Riesz that such measures are absolutely continuous, that is have the form 
$d\nu(\theta)=g(e^{i\theta})\,d\theta$ 
with $g\in H^1$. The Hahn-Banach theorem 
implies that $H^1$ is dual {\it via} (\ref{pairdual}) to the quotient space
$C(\TT)/\mathcal{A}_0$ \cite[chap. IV, sec. 1]{Garnett}. 
Equivalently, ${\bar H}^1_0$ is dual to $C(\TT)/\overline{\mathcal{A}}$ under the pairing
arising from the line integral~:
\begin{equation}
\label{pairbar}
(\dot{f},F)=
 \frac{1}{2i\pi}\int_{\TT} f(\xi)F(\xi)\,d\xi\,,
\end{equation} 
where $F$ belongs to ${\bar H}^1_0$ and $\dot{f}$
indicates the equivalence class of $f\in C(\TT)$
modulo $\overline{\mathcal{A}}$. 
Therefore, contrary to $L^1(\TT)$, the spaces
$H^1$ and ${\bar H}^1_0$ enjoy a weak-* compactness property
of their unit ball.

We define the analytic and anti-analytic projections
${\bf P}_+$ and ${\bf P}_-$ on Fourier series by~:
\[
{\bf P}_+\left(\sum_{n=-\infty}^{\infty}a_ne^{in\theta}\right)=
\sum_{n=0}^{\infty}a_ne^{in\theta},~~~~
{\bf P}_-\left(\sum_{n=-\infty}^{\infty}a_ne^{in\theta}\right)=
\sum_{n=-\infty}^{-1}a_ne^{in\theta}\,.
\]
It is a theorem of M. Riesz theorem \cite[ch. III, sec, 1]{Garnett}
 that ${\bf P}_+:L^p\to H^p$ and ${\bf P}_-:L^p\to \bar{H}_0^p$
are bounded for $1<p<\infty$, in which
case they coincide with the Cauchy projections:
\begin{equation}
\label{Cauchyproj}
{\bf P}_+(h)(z)=\frac{1}{2i\pi}\int_{\TT}\frac{h(\xi)}{\xi-z}\,d\xi,
~~~z\in\DD,~~~~~~~
{\bf P}_-(h)(s)=\frac{1}{2i\pi}\int_{\TT}\frac{h(\xi)}{s-\xi}\,d\xi,
~~~~s\in\hat{\CC}\setminus\overline{\DD}.
\end{equation}
When restricted to $L^2(\TT)$, the projections ${\bf P}_+$ and
${\bf P}_-$ are just the
orthogonal projections onto $H^2$ and $\bar{H}_0^2$ respectively. 
Although ${\bf P}_\pm(h)$ needs not be the Fourier series of a function
when $h$ is merely in $L^1(\TT)$, it is Abel summable 
almost everywhere to a function lying in $L^s(\TT)$
for $0<s<1$ and it can
be interpreted as a function in the
Hardy space of exponent $s$ that we did not introduce 
\cite[cor. to thm 3.2]{Duren}. To us it will be sufficient, 
when $h\in L^1$, to regard  ${\bf P}_\pm(f)$ as the Fourier 
series of a distribution. 

Finally, we let $P_n$ denote throughout the space of complex algebraic polynomials of degree at most $n$. Clearly,
$P_n\subset H^p$ for all $p$.


\section{Two extremal problems}
\label{wposed}
We first state the polynomial problem discussed in Section  \ref{sec1}. 
We call it $PBEP(n)$ for ``Polynomial Bounded Extremal Problem'':
\begin{description}
\item[PBEP(n)] \mbox{}\\
For $f \in L^2(I)$, find $k_n \in P_n$ such that
$|k_n(e^{i \theta})| \leq 1$ for a.e. $e^{i \theta} \in J$
and 
\begin{equation} 
\|f - k_n\|_{L^2(I)} = \inf_{\stackrel{g \in P_n}{|g| \leq 1 \,
    \mbox{a.e. on} \, J}} \|f - g\|_{L^2(I)} \, .
\label{eq-approxp}
\end{equation}  
\end{description}
Next, we state the analytic problem from Section \ref{sec1} that we call
$ABEP$ for ``Analytic Bounded Extremal Problem'':
\begin{description}
\item[ABEP] \mbox{}\\
Given $f \in L^2(I)$, find $g_0 \in H^2$ such that 
$|g_0(e^{i \theta})| \leq 1$ a.e. on $J$ and
\begin{equation} 
\|f - g_0\|_{L^2(I)} = \inf_{\stackrel{g \in H^2}{|g| \leq 1 \,
    \mbox{a.e. on} \, J}} \|f - g\|_{L^2(I)} \, .
\label{eq-approx}
\end{equation}  
\end{description}
Note that, in $ABEP$, the constraint 
$|g| \leq 1$ on $J$ could be replaced by $|g| \leq \rho$ where $\rho$ is a positive function in
$L^2(J)$. For if $\log \rho\in L^1(J)$ then, denoting
by $w_{1\vee(1/\rho)}$ the outer factor having modulus $1$ on $I$ and 
$1/\rho$ on $J$, we find that
$g\in H^2$ satisfies $|g |\leq\rho$ on $J$ if and only if
$h=gw_{1\vee(1/\rho)}$ lies in $H^2$ and satisfies $|h|\leq1$ on $J$.
It is so because, for $g$ as indicated, 
$h$ lies in the Nevanlinna class by construction and
$|h|_{|_I}=|g|_{|_I}$ while $|h|_{|_J}=|g|_{|_J}/\rho$.
If, however, $\log \rho\notin L^1(J)$, then we must have
$\int_J\log\rho=-\infty$ because $\rho\in L^2(J)$, consequently 
the set of candidate approximants 
reduces to $\{0\}$ anyway because a nonzero Hardy function has summable
$\log$-modulus. Altogether, it is thus equivalent  
to consider
$ABEP$ for the product $f$ times ${(w_{1\vee\rho^{-1}})}_{|_I}$.
A similar argument shows that we could replace  the error criterion
$\|.\|_{L^2(I)}$
by a weighted norm $\|.\|_{L^2(I,w)}$ for some weight  $w$ which
is non-negative and invertible in $L^\infty(I)$. Then, the problem reduces to $ABEP$ for
$f(w_{\rho^{1/2}\vee0})_{_I}$. 

Such equivalences do not hold for $PBEP(n)$ because the polynomial character of 
$k_n$ is not preserved under multiplication by outer factors. Still, 
the results to come continue to hold if we replace in $PBEP(n)$ the constraint
$|k_n| \leq 1$ by $|k_n| \leq \rho$ on $J$ and the criterion 
$\|.\|_{L^2(I)}$ by $\|.\|_{L^2(I,w)}$, provided that $\rho\in C(J)$ and
that $w$ is invertible in  $L^\infty(I)$. Indeed, we leave it to the reader to check
that proofs go through with obvious modifications.

After these preliminaries, we are ready to state a basic existence and uniqueness result.

\begin{theorem}.
\label{exunc}
Problems $PBEP(n)$ and $ABEP$ have a unique solution. Moreover, the solution
$g_0$ to $ABEP$ satisfies 
$|g_0|=1$ almost everywhere on $J$, unless 
$f=g_{|_I}$ for some $g\in H^2$ such that 
$\|g\|_{L^\infty(J)}\leq 1$.
\end{theorem}
{\sl Proof.}
Consider the sets
\[
E_n=\{g_{|_I}:\ g\in P_n,\  \|g\|_{L^\infty(J)}\leq1\},
\]
\[
F=\{g_{|_I}:\ g\in H^2,\  \|g\|_{L^\infty(J)}\leq1\}.
\]
Clearly $E_n\subset F$  are convex and nonempty subsets of
$L^2(I)$, as they contain 0.
To prove existence and uniqueness, it is therefore enough to show they are 
closed, for we can appeal then to well-known properties of the
projection on a closed convex set in a Hilbert space.
Since $E_n=P_n\cap F$, it is enough in fact to show that $F$ is closed.
For this, let $g_m$ be a sequence in $H^2$ with $|g_m|_{|_J}\leq 1$ and such 
that $(g_m)_{|_I}$ converges in $L^2(I)$. Obviously $g_m$ is a  bounded 
sequence in 
$L^2(\TT)$, some subsequence of which converges weakly to $h\in H^2$. We continue to denote this subsequence with $g_m$.
The restrictions $(g_m)_{|_I}$ a fortiori converge weakly to
$h_{|_I}$ in $L^2(I)$, and since the strong and the weak limit must coincide 
when both exist we find that $(g_m)_{|_I}$ converges to $h_{|_I}$ in $L^2(I)$.
Besides,  $(g_m)_{|_J}$ is contained in the unit ball of $L^\infty(J)$
which is dual to $L^1(J)$, hence some subsequence (again denoted by 
 $(g_m)_{|_J}$)
converges weak-* to some $h_1\in L^\infty(J)$ with
$\|h_1\|_{L^\infty(J)}\leq1$. But since $(g_m)_{|_J}$ also converges weakly to
$h_{|_J}$ in $L^2(J)$, we have that
\[\langle h_1,\varphi\rangle_{J}=
\lim_{m\to\infty}\langle g_m,\varphi\rangle_{J}=\langle h_{|_J},
\varphi\rangle_{J}
\]
for all $\varphi\in L^2(J)$ which is dense in $L^1(J)$. Consequently
$h_1=h_{|_J}$, thereby showing that $\|h\|_{L^\infty(J)}\leq1$, 
which proves that $F$ is closed.

Assume now that $f$ is not the trace on $I$ of an $H^2$-function which is less than 1 in modulus on $I$. To prove that $|g_0|=1$ a.e. on $J$, 
we argue by contradiction.
If not, there is a compact set $K$ of positive measure,
lying interior to $J$, such that
$\|g_0\|_{L^\infty(K)}\leq 1-\delta$ for some $0<\delta<1$;
it is so because, by hypothesis, $J$ must consist of finitely many 
closed arcs, of which one at least has nonempty interior. 
For $K'$ an arbitrary subset of $K$,
consider the Riesz-Herglotz transform of its
characteristic function:
\begin{equation}
h_{K'}(z) = \frac{1}{2 \, \pi} \, \int_{K'} \frac{e^{i\theta} + 
z}{e^{i\theta} - z} \, \, d \theta \, 
,~~~~z\in\DD,
\label{defh}
\end{equation}
and put $w_{t} = \exp (th_{K'})$ for $t \in \RR$,
which is the outer function with modulus 
$\exp t$ on $K'$ and  $1$ elsewhere. By construction, 
$g_0 w_t$ 
is a candidate approximant in
$ABEP$ for all $t<-\log(1-\delta)$. Thus, the map
$t\mapsto \|f-g_0w_t\|^2_{L^2(I)}$ 
attains a minimum at $t=0$. Because $K$ is at strictly 
positive distance from $I$, we may
differentiate this expression with respect to $t$ under the integral sign
and equate the derivative at $t=0$ to zero which gives us
$2\mbox{\rm Re} \langle f-g_0,h_{K'} g_0 \rangle_I=0$. Replacing $g_0w_{t}$ by $ig_0w_{t}$,
which is a candidate approximant as well, we get a similar equation
for the imaginary part so that
\begin{equation}
\label{OSu}
0=\langle f-g_0\,,\,h_{K'} g_0\rangle_I=\langle (f-g_0)\bar{g}_0\,,\,
h_{K'}\rangle_I.
\end{equation}
Let $e^{it_0}$ be a density point of $K$ and $I_l$ 
the arc centered at $e^{it_0}$ of length $l$, so that
$|I_l\cap K)|/l\to1$ as $l\to0$. 
Since
\begin{equation}
\label{comparfrac}
\left|\frac{e^{it}+e^{i\theta}}{e^{it}-e^{i\theta}}-
\frac{e^{it_0}+e^{i\theta}}{e^{it_0}-e^{i\theta}}\right|
\leq \frac{2l}{\mbox{\rm dist}^2(K,I)}~~~~\mathrm{for}~e^{it}\in I_l\cap K,~~e^{i\theta}
\in I,
\end{equation}
it follows by dominated convergence that 
\[
\lim_{l\to0}\frac{1}{|I_l\cap K|}\int_{I_l\cap K}
\left|\frac{e^{it}+e^{i\theta}}{e^{it}-e^{i\theta}}-
\frac{e^{it_0}+e^{i\theta}}{e^{it_0}-e^{i\theta}}\right|dt=0,\quad
\mbox{\rm uniformly w.r. to\ }e^{i\theta}\in I,
\]
and therefore that
\[
\lim_{l\to0}
\frac{h_{I_l\cap K}(e^{i\theta})}{|I_l\cap K|}=
\frac{e^{it_0}+e^{i\theta}}{e^{it_0}-e^{i\theta}}
\quad
\mbox{\rm uniformly w.r. to\ }e^{i\theta}\in I.\]
Applying now \eqref{OSu} with $K'=I_l\cap K$ and taking into account that
$(e^{it_0}+e^{i\theta})/(e^{it_0}-e^{i\theta})$ is pure imaginary on $I$,
we find in the limit, as $l\to0$ that 
\begin{equation}
\label{RHv}
\frac{1}{2\pi}\int_I \frac{e^{it_0}+e^{i\theta}}{e^{it_0}-e^{i\theta}}
\Bigl((f-g_0)
\bar{g}_0\Bigr)(e^{i\theta})\,d\theta=0.
\end{equation}
Next, let us consider the function
\[
F(z)=\frac{1}{2\pi}\int_I \frac{e^{i\theta}+z}{e^{i\theta}-z}
\Bigl((f-g_0)
\bar{g}_0\Bigr)(e^{i\theta})\,d\theta=
-\frac{1}{2\pi}\int_I\Bigl((f-g_0)
\bar{g}_0\Bigr)(e^{i\theta})d\theta+
\frac{1}{i\pi}\int_I \frac{\Bigl((f-g_0)
\bar{g}_0\Bigr)(\xi)d\xi}{\xi-z}
\]
which is the sum of a 
constant and of twice the Cauchy integral of
$(f-(g_0)_{|_I})(\bar{g}_0)_{|_I}\in L^1(I)$, hence is 
analytic in $\hat{\CC}\setminus I$. Equation \eqref{RHv} means
that $F$ vanishes at every density point of $K$, and since a.e. point in 
$K$ is a density point  $F$ must vanish identically
because its zeros accumulate in the interior of $J$.
Denoting by $F^+$ and $F^-$ the nontangential limits of $F$ from sequences of points in $\DD$ or $C \setminus \overline{\DD}$ respectively,
we now get from  the 
Plemelj-Sokhotski formulas \cite[ch. III]{Garnett} that
\[0=F^+(\xi)-F^-(\xi)=(f-g_0)(\xi)\overline{g_0(\xi)},\qquad
a.e. \ \xi\in I.
\]
Thus, either $g_0$ is nonzero a.e. on $I$, in 
which case $f=(g_0)_{|_I}$ and we reach the desired contradiction, or else
$g_0\equiv0$. In the latter case, if we put $\mbox{\rm id}$ for the identity 
map on $\TT$, we find that
$t\mapsto \|f-t\,\mbox{\rm id}^k\|^2_{L^2(I)}$ has a minimum at $t=0$
for each integer $k\geq0$,
since $e^{i\theta}\mapsto t e^{ik\theta}$ is a candidate approximant
for $t\in[-1,1]$. Differentiating with respect to $t$ and expressing that the derivative at $t=0$ is zero, we deduce that all Fourier coefficients of 
non-negative index of  $(f-(g_0)_{|_I})\vee0$ do vanish.
This means this last function lies in ${\bar H}^2$, but 
as it vanishes on $J$ it is identically zero, therefore  $f=(g_0)_{|_I}$ in all cases.
\hfill \boite\\
{\bf Remark:} the theorem shows that the constraint $|g_0|\leq1$ on $J$ 
is saturated in a very strong sense for problem $ABEP$, namely
$|g_0|=1$ a.e. on $J$ unless $f$ is
already the trace of the solution on $I$.
In contrast, it is not true that $\|k_n\|_{L^\infty(J)}=1$ unless 
$f=g_{|_I}$ for some $g\in P_n$ such that 
$\|g\|_{L^\infty(J)}<1$. To see this, observe that the set $E_n$ is not only 
closed but compact. Indeed, if we pick distinct points 
$\xi_1,\cdots,\xi_{n+1}$ in 
$J$ and form the Lagrange interpolation polynomials 
$L_j\in P_n$ such that 
$L_j(\xi_j)=1$ and $L_j(\xi_\ell)=0$ if $\ell\neq j$, we get
a basis of $P_n$
in which the coordinates of every $g\in P_n$ meeting 
$\|g\|_{L^\infty(J)}\leq1$ are bounded by 1 in modulus.
Hence $E_n$ is bounded in $P_n$, and since it is closed by the proof
of Theorem \ref{exunc} it is compact.
Thus, each $f\in L^2(I)$ has a best approximant from $E_n$, 
and if $(p_n)_{_I}$ is a best 
approximant to $f$ with $p_n\in P_n$, then
for $\lambda>\|p_n\|_{L^\infty(J)}$ we find that $p_n/\lambda$ is a best
approximant to $f/\lambda$ in $L^2(I)$ which is strictly less than 1 on $J$.
This justifies the remark.

\section{Critical point equations  and convergence of approximants}
\label{critpsec}
At this point,  it is worth recalling informally some
basic principles  from convex optimization, for which the reader may consult
\cite{BorLew}. The solution  to a strictly  convex minimization problem
is characterized by
a variational inequality expressing that the
{\it criterium } increases  under admissible increments of the variable.
If the problem is smooth enough, such increments admit  a tangent space at the point under consideration ({\it i.e.} the solution) in the variable space. 
We term  it the tangent space to the constraints,
and its orthogonal in the dual space to the variable space is called
the orthogonal space to the constraints (at the point under consideration).
The variation of the objective function
must vanish on the tangent space to the constraints 
to the first order, thereby giving rise
to the so-called  \emph{critical point equation}.
It says that the gradient of the objective function,
viewed as an element of the dual space to the variable space, lies in the 
orthogonal space to the constraints.
If a basis of the latter is chosen,
the coordinates of the gradient in this basis 
are known as the  \emph{Lagrange parameters}. 
More generally, one can form the Lagrangian which is a function of the 
variable and of the 
Lagrange parameters, not necessarily optimal ones. It is obtained 
by adding the gradient of the criterion,  at the considered value
of the variable,
with the member of the orthogonal space to the constraints 
defined by the chosen Lagrange parameters.
By what precedes, the Lagrangian must vanish at the solution 
for appropriate values of the
Lagrange parameters. One can further define a function of the 
Lagrange parameters 
only, by minimizing the Lagrangian with respect to the variable.
This results in a concave function which gets maximized at the optimal value 
of  the Lagrange parameters for the original problem. 
This way, one reduces the original 
constrained convex  minimization problem to an unsconstrained concave
maximization problem, called the dual problem. In an infinite-dimensional
context, the arguments needed to put this program to work may be quite subtle.

Below we derive the critical point equation for  $PBEP(n)$ described in
\eqref{eq-approxp}.
For $g\in P_n$ define $$E(g)=\{x \in
J,\,\,|g(x)|=||g||_{L^{\infty}(J)}\},$$ which is the set of
extremal  points of $g$ on $J$.

\begin{theorem}
\label{th_dimf}
A polynomial $g\in P_n$ is the solution to $PBEP(n)$  iff the following 
two conditions hold:
\begin{itemize}
\item $||g||_{L^{\infty}(J)} \leq 1$, \\
\item there exists
a set of $r$ distinct points $x_1,
\cdots,x_r\in E(g)$  and 
non-negative real numbers $\lambda_1,\cdots,\lambda_r$,
with $0\leq r\leq 2n+2$, such that
\begin{equation}
 \langle g-f,h \rangle_{I}+\sum_{j=1}^{r} \lambda_j
g(x_j)\overline{ h(x_j)} =0,\qquad \forall h \in P_n. 
\label{carac_fbep}
\end{equation}
\end{itemize}
Moreover the $\lambda_j$'s meet  the following bound 
\begin{equation}
\label{lambda_bound}
\sum_{j=1}^{r} \lambda_j \leq 2||f||_{L^2(I)}^2.
\end{equation}
We emphasize that the set of extremal points $\{x_j,\,\,j=1,\dots ,r\}$ 
is possibly empty (i.e $r=0$). 
\end{theorem}
{\sl Proof.}
Suppose $g$ verifies the two conditions and differs from the solution
$k_n$. Set
$h=k_n-g\in P_n$ and observe that 
\begin{equation}
\label{neg} 
\mbox{\rm Re}\left(g(x_i)\overline{h(x_i)}\right)=\mbox{\rm Re}\left(g(x_i)\overline{ k_n(x_i)}-1\right) \leq 0,
\,\,\,i=1\dots r .
\end{equation} 
From the uniqueness and optimality of $k_n$ we deduce that
\begin{equation}
\nonumber
\begin{split} 
||k_n-f||^2_{L^2(I)}&=||g-f+h||^2_{L^2(I)} \\
                  &=||g-f||^2_{L^2(I)}+||h||^2_{L^2(I)}+2\mbox{\rm Re} \langle g-f,h \rangle_{I} \\
                  &<||g-f||_{L^2(I)}^2.
\end{split}
\end{equation}
Consequently $ \mbox{\rm Re} \langle g-f,h \rangle_{L^2(I)}<0$ which, combined with
(\ref{neg}), contradicts (\ref{carac_fbep}).

Conversely, suppose that $g$ is the solution
to $PBEP(n)$ and let $\phi_0$ be the
$\Bbb R$-linear forms on $P_n$ given by
$$
 \,\,\phi_0(h)=\mbox{\rm Re}\langle g-f,h\rangle_{I},\qquad  h \in P_n.
$$
For each extremal point $x\in E(g)$, define further a $\Bbb R$-linear
form $\phi_x$ by
$$\,\, \phi_x(h)=\mbox{\rm Re}\left(g(x)\overline{h(x)}\right),\quad
h\in P_n.$$
Put $K$ for the union of these forms:
$$K=\{\phi_0\} \cup \{\phi_x,\,\,x \in E(g)\}.$$ 
If we let  $P_n^{\RR}$  indicate $P_n$ viewed as
a real vector space, $K$ is a subset of  the dual $(P_n^{\RR})^*$.
As $J$ is closed by definition,
simple inspection shows that $K$ is closed
and bounded in $(P_n^{\RR})^*$ 
(it is in fact finite unless $g$ is a constant), hence 
it is compact and so is its convex hull
$\hat{K}$ as $(P_n^\RR)^*$ is finite-dimensional.
Suppose for a contradiction that $0 \not \in \hat{K}$. Then,
since $(P_n^{\RR})^{**}=P_n^{\RR}$  because
$P_n^\RR$ is finite-dimensional,
there exists by the Hahn-Banach theorem an $h_0 \in P_n$ such that,
$$ \,\, \phi(h_0) \geq \tau > 0, \quad \forall \phi \in \hat K . $$  
The latter and the continuity of $g$ and $h_0$ ensure the existence of
a neighborhood $V$ of $E(g)$ on $\TT$ 
such that for $x$ in  $U=J \cap V$ we have 
$\mbox{\rm Re} \left(g(x)\overline{h_0(x)}\right) \geq \frac{\tau}{2}>0$,
whereas 
for $x$ in $J\backslash U$ it holds that  $|g(x)| \leq 1-\delta
\,\,$ for some $\delta>0$.
Clearly, for $\epsilon >0 $ with $\epsilon
||h_0||_{L^{\infty}(J)}< \delta$, we get  that
\begin{equation}
\sup_{J \backslash U} |g(x) - \epsilon h_0(x)| \leq 1. \label{surcU} 
\end{equation}
Moreover, assuming without loss of generality that $\epsilon<1$, 
it holds for  $x \in U$ that
\begin{equation}
\nonumber
\begin{split}
|g(x)-\epsilon h_0(x)|^2 &=
|g(x)|^2-2\mbox{\rm Re} \left(g(x)\overline{h_0(x)}\right)+\epsilon^2|h_0(x)|^2 \\
&\leq |g(x)|^2-2\mbox{\rm Re} \left(\epsilon g(x)\overline{h_0(x)} \right)+\epsilon^2|h_0(x)|^2 \\
&\leq 1-\epsilon \tau+\epsilon^2||h_0||_{L^{\infty}(J)}^2. \\
\end{split}
\end{equation}
The latter combined with (\ref{surcU}) shows that, 
for $\epsilon$ sufficiently
small, we have
\begin{equation}
||g-\epsilon h_0||_{L^{\infty}(J)}\leq 1 \label{surJ}.
\end{equation}
However, since
\begin{equation}
\begin{split}
||f-g-\epsilon h_0||_{L^2(J)}^2&=||f-g||_{L^2(J)}^2-2\epsilon
\phi_0(h_0)+\epsilon^2||h_0||_{L^2(J)}^2 \\
&\leq ||f-g||_{L^2(J)}^2-2\epsilon \tau + \epsilon^2||h_0||_{L^2(J)}^2 ,
\end{split}
\end{equation}
we deduce in view of  (\ref{surJ}) that for $\epsilon$ small enough
the polynomial $g-\epsilon h_0$ performs better than $g$ 
in $FBEP$, thereby contradicting optimality. Hence $0 \in \hat K$,
therefore by Carathedory's theorem \cite[ch. 1, sec. 5]{Cheney}
there are $r'$
elements $\gamma_j$ of $K$, with $1\leq r' \leq 2(n+1)+1$
(the real dimension of $P_n^\RR$ plus one), such that
\begin{equation}
\label{convex_comb}
\sum_{j=1}^{r'} \alpha_j \gamma_j =0 \\
\end{equation}
for some positive  $\alpha_j$ satisfying $\sum
\alpha_j=1$.
Of necessity $\phi_0$ is a $\gamma_j$, otherwise
evaluating (\ref{convex_comb}) at $g$ yields the absurd conclusion that
$$ 0= \sum_{j=1}^{r'} \alpha_j \gamma_j(g)= \sum_{j=1}^{r'} \alpha_j
|g(x_j)|^2=1.$$ 
Equation (\ref{convex_comb}) can therefore be rewritten
as
$$\,\ \alpha_1 \mbox{\rm Re} \langle f-g,h \rangle_{I}+\sum_{j=2}^{r'} \alpha_j
\mbox{\rm Re}(g(x_j)\overline{h(x_j)}) = 0\quad
\forall h \in P_n,\quad \alpha_1\neq0 .$$ Dividing 
by $\alpha_1$ and noting that the last equation is also true
with $ih$ instead of $h$ yields 
(\ref{carac_fbep}) with $r=r'-1$.
Finally, replacing $h$ by $g$ in 
(\ref{carac_fbep}) we obtain

\begin{equation}
\nonumber
\begin{split}
\sum_{j=1}^{r}|\lambda_j|=\sum_{j=1}^{r}\lambda_j|=\langle f-g,g\rangle_{I}  
          &\leq \,\, \langle f-g,f-g\rangle_{I}+|\langle f-g,f\rangle_{I}|\\
          &\leq ||f-g||_{L^2(I)}^2+||f-g||_{L^2(I)}||f||_{L^2(I)} \\
          &\leq 2 ||f||_{L^2(I)}^2 
\end{split}
\end{equation}

where the next to last majorization uses the Schwarz inequality and the last
that $0$ is a candidate approximant for $PBEP(n)$ whereas  $g$ is the optimum. 
\hfill \boite\\
 
The next result describes the behavior of $k_n$ when
$n$ goes to infinity, in connection with the solution $g_0$ to $ABEP$.

\begin{theorem}
\label{conv_n}
Let $k_n$ be the solution to $PBEP(n)$ defined in \eqref{eq-approxp},
 and $g_0$ the solution to
$ABEP$ described in \eqref{eq-approx}.
When $n\to\infty$, the sequence $(k_n)_{|_I}$ converges to $(g_0)_{|_I}$ in
$L^2(I)$ , and  the sequence $(k_n)_{|_J}$ 
converges to $(g_0)_{|_J}$ in
the weak-* topology of $L^{\infty}(J)$, as well as in
$L^p(J)$-norm for $1 \leq p < \infty$ if $f$ is not the race on $I$
of a $H^2$-function which is at most 1 in modulus on $J$. 
Altogether this amounts to:
\begin{equation}
\label{convI}
\lim_{n \rightarrow \infty} ||g_0 - k_n||_{L^p(\TT)}=0,\quad 1\leq p\leq2, 
\end{equation}
\begin{equation}
\label{convJ}
\,\,\lim_{n \rightarrow \infty} \langle k_n,h \rangle_{J}=\langle g_0,h \rangle_{J}\quad \forall h \in L^1(J),
\end{equation}
\begin{equation}
\label{convpJ}
\mbox{\rm if $f\neq g_0$ on $I$},\ \lim_{n \rightarrow \infty} ||g_0 - k_n||_{L^p(J)}=0, \quad
1 \leq p < \infty. 
\end{equation}
\end{theorem}
{\sl Proof.}
Our first objective is to show that $g_0$ can be approximated 
arbitrary close 
in $L^2(I)$ by
polynomials that remain bounded by 1 in modulus on $J$.
 By hypothesis $I$ is the finite union of
$N\geq1$ open disjoint sub-arcs of $\TT$. Without loss of generality,
it can thus  be written as
$$I=\bigcup_{i=1}^{N} (e^{ia_i},e^{ib_i}),
\qquad0=a_1 \leq b_1 \leq a_2 \dots \leq b_N \leq 2\pi.$$
Let $(\epsilon_n)$ be a sequence of positive real numbers
decreasing to  0.
We define  a sequence $(v_n)$ in $H^2$ by
\begin{equation}
\nonumber
\begin{split}
v_n(z)&=
g_0(z) \exp\left(-\frac{1}{2\pi}\left(\sum_{i=1}^{N}\int_{a_i}^{a_i+\epsilon_n}\frac{e^{it}+z}{e^{it}-z}\log|g_0|dt+\int_{b_i-\epsilon_n}^{b_i}\frac{e^{it}+z}{e^{it}-z}\log|g_0|dt\right)\right)
\label{varg0}
\end{split}
\end{equation}
Note that indeed $v_n\in H^2$ for $n$ large enough because
then it has the same modulus as $g_0$ except over the arcs 
$(a_i,a_i+\epsilon_n)$
and $(b_i-\epsilon_n,b_i)$ where it has modulus 1.
\emph{We claim} that $(v_n)_{|_I}$ converges to $g_0$ in $L^2(I)$ as 
$n\to\infty$.
To see this, observe that $v_n$ converges a.e. on 
$I$ to $g_0$, for  each  $z\in I$
remains at some distance from 
the sub-arcs
$(a_i,a_i+\epsilon_n)$ and $(b_i,b_i+\epsilon_n)$ 
for all $n$ sufficiently
large, hence the argument of the exponential in \eqref{varg0} converges to zero
as $n\to\infty$ by absolute continuity of $\log|g_0|dt$. 
Now, we remark that by
construction $|v_n| \leq |g_0|+1$, 
hence by dominated convergence, we get that
$$\lim_{n \rightarrow \infty} ||g_0-v_n||_{L^2(I)}=0.$$      
\emph{This proves the claim}.
Now, let $\epsilon >0$ and $0<\alpha<1$ such that $||g_0-\alpha
g_0||_{L^2(I)}\leq \frac{\epsilon}{4}$. Let also
$n_0$ be so large 
that $||v_{n_0}-g_0||_{L^2(I)}\leq \frac{\epsilon}{4}$. 
For $0<r<1$
define $u_r\in\mathcal{A}$ (the disk algebra) by $u_r(z)=v_{n_0}(rz)$ so that, by Poisson representation,
$$u_r(e^{i\theta})=\int_{\TT}P_r(\theta-t)v_{n_0}(re^{it})dt,$$ 
where $P_r$ is the Poisson kernel. 
Whenever $e^{i\phi} \in J$, we note by construction that
$|v_n|=1$ a.e on the sub-arc
$(e^{i(\phi-\epsilon_{n_0})},e^{i(\phi+\epsilon_{n_0}}))$. This is to the effect 
    that 
\begin{equation}
\nonumber
\begin{split}
|u_r(e^{i\phi})| & \leq  \int_{\TT}P_r(\phi-t)|v_{n_0}(re^{it})|dt \\
           & \leq  P_r(\epsilon_{n_0})\int_{\TT}|v_{n_0}(re^{it})|dt+
\int_{-\epsilon_{n_0}}^{+\epsilon_{n_0}} P_r(t)dt\\
           & \leq P_r(\epsilon_{n_0})||v_{n_0}||_{L^1(\TT)}+
          1 \leq P_r(\epsilon_{n_0})||v_{n_0}||_{L^2(\TT)}+
          1 
\end{split}
\end{equation}
by H\"older's inequality.
Hence, for $r$ sufficiently close to $1$, we certainly have that
$|u_r| \leq 1/\alpha^2$ on
$J$ and otherwise that $||u_r-v_{n_0}||^2_{L^2(I)} \leq \frac{\epsilon}{4}$
since $u_r\to v_{n_0}$ in $H^2$. Finally,
call $q$ the truncated Taylor expansion of $u_r$ (which converges
uniformly to the latter on $\TT$), 
where the order of truncation has been chosen large
enough  to ensure that $|q| \leq  1/\alpha$ on $J$ and that 
$||q-u_r||^2_{L^2(I)} \leq \frac{\epsilon}{4}$. Then, we have that
\begin{equation}
\nonumber
\begin{split}
||\alpha q - g_0||_{L^2(I)} &\leq 
\alpha \left(||q-u_r||_{L^2(I)}+||u_r-v_{n_0}||_{L^2(I)}+
||v_{n_0}-g_0||_{L^2(I)}\right)+||g_0-\alpha g_0||_{L^2(I)} \\
                            &\leq  \epsilon .\\
\end{split}
\end{equation}
Thus, we have found a polynomial (namely $\alpha q$) which is bounded by 1 
in modulus on $J$ and close by $\epsilon$ to $g_0$ in $L^2(I)$.
By comparison, this immediately implies that
\begin{equation}
\label{eps4}
\lim_{n \rightarrow \infty}||f-k_n||_{L^2(I)}=||f-g_0||_{L^2(I)},
\end{equation} 
from which \eqref{convI} follows by H\"older's inequality. Moreover,
being bounded in $H^2$, the sequence $(k_n)$ has  a weakly
convergent sub-sequence. The traces on $J$ of this subsequence are
in fact bounded by 1 in $L^{\infty}(J)$-norm, hence up to another
subsequence we obtain $(k_{n_m})$ converging also
in the 
weak-* sense on $J$. Let $g$ be the weak limit ($H^2$ sense) of 
$k_{n_m}$, and observe that $g_{|_J}$ is necessarily 
the weak-* limit of $(k_{n_m})_{|_J}$ in $L^\infty(J)$, as follows by integrating 
against functions from $L^2(J)$ which is dense in $L^1(J)$.
Since  balls
are weak-* closed in $L^\infty(J)$, we have that $||g||_{L^{\infty}(J)}
\leq 1$, and it follows from (\ref{eps4}) that
$||f-g||_{L^2(I)}=||f-g_0||_{L^2(I)}$. Thus,
$g=g_0$ by the uniqueness part of Theorem \ref{exunc}.
Finally,
if $f\neq g_0$ on $J$, then we know from Theorem \ref{exunc} that
$|g_0|=1$ a.e. on $J$.
In this case,
(\ref{eps4}) implies
that $\limsup ||k_{n_m}||_{L^2(\TT)} \leq ||g_0||_{L^2(\TT)}$, and since the 
norm of the weak limit is no less than the limit of the norms it follows that
$(k_{n_m})_{|_J}$ converges strongly to $(g_0)_{|_J}$ in the 
strictly convex space $L^2(J)$. The same 
reasoning applies in $L^p(J)$ for $1<p<\infty$.
Finally we remark
that the preceding arguments hold true when $k_n$ is replaced by
any  subsequence of itself; hence $k_n$ contains no
subsequence not converging to $g_0$ in the sense stated before,
which achieves the proof.
\hfill \boite\\ 

We come now to an analog of theorem \ref{th_dimf} in the infinite
dimensional case. We define $H^{2,\infty}_J$ and $H^{2,1}_I$ to be the
following vector spaces:
$$ H^{2,\infty}_J = \{h \in H^2,\,\,||h||_{L^{\infty}(J)}<\infty\},$$
$$ H^{2,1}_I=\{h \in H^1,\,\,||h||_{L^2(I)}<\infty\},$$
endowed with the natural norms.
We begin with an elementary lemma.

\begin{lemma}
\label{lemme_integr}
Let $v \in L^{1}(J)$ such that ${\bf P}_+(0 \vee v) \in H^{2,1}_I$. Then: 
$$ \forall h \in H^{2,\infty}_J,\,\,\langle{\bf P}_+(
0 \vee v),h \rangle_{\TT}= \langle v,h \rangle_{J}. $$
\end{lemma}
{\sl Proof.}
Let $u$ be the function defined on $\TT$ by $$u=(0 \wedge v) -{\bf P}_+(
0 \vee v).$$ By assumption $u\in L^1(\TT)$, and
by its very definition all Fourier coefficients of $u$
of non-negative index vanish. Hence $u\in \bar{H}_0^1$, and 
since it is $L^2$ integrable on $I$ where it coincides with 
$-{\bf P}_+(0\vee v)$, 
.we conclude that $\overline{u} \in
H^{2,1}_I$ and that $\overline{u}(0)=0$ 
Now, for  $h \in H^{2,\infty}_J$ we have that
\begin{equation}
\begin{split}
\langle v\chi_J,h \rangle_{\TT}&=\,\,\langle u,h \rangle_{\TT}+ \langle {\bf P}_+(0 \vee v ),h \rangle_{\TT} \\
                      &=\overline{u}(0)h(0)+\langle {\bf P}_+(0 \vee v ),h \rangle_{\TT} \\
                      &=\,\,\langle{\bf P}_+(0 \vee v),h\rangle_{\TT} \\
\end{split}
\end{equation}

where the second equality follows from the Cauchy formula
because $(\overline{u}h) \in H^1$.
\hfill \boite\\

\begin{theorem}
\label{th_dimi}
Suppose that $f \in L^2(I)$ is not the trace on $I$ of a $H^2$-function 
of modulus less or equal to $1$ a.e on $J$. Then, $g\in H^2$ is the 
solution to $ABEP$ iff the following 
two conditions hold.
\begin{itemize}
\item \,\, $|g(e^{i\theta})|=1$\,\, 
for a.e. $ e^{i\theta} \in J$,\\
\item there exists
a nonnegative real function $\lambda \in L^1_\RR(J)$ such that, 
\begin{equation}
\forall h \in H^{2,\infty}_J,\, \langle g-f,h \rangle_{I}+\langle \lambda g,h \rangle_{J}=0. 
\label{carac_bep}
\end{equation}
\end{itemize} 
\end{theorem}
{\sl Proof.}
Suppose $g$ verifies the two conditions and differs from
$g_0$. Set $h=(g_0-g) \in H^{2,\infty}_J $ and observe that
\begin{equation}
\label{neg1}
\mbox{\rm Re} \langle \lambda g,h \rangle_{J}=\frac{1}{2\pi}\int_{J} \lambda(\mbox{\rm Re}(\overline{g}g_0)-1)
\leq 0. 
\end{equation}
In another connection, since $-h$ is an admissible increment from $g_0$,
the variational inequality characterizing
the projection onto a closed convex set gives us ({\it cf.} Theorem \ref{exunc})
$\mbox{\rm Re}\langle g_0-f,h \rangle_{I}\leq0$, whence
$$ \mbox{\rm Re} \langle g-f,h \rangle_{I}=
\mbox{\rm Re} \langle g_0-f,h \rangle_{I}-\langle h,h \rangle_I
<0$$ which, combined with
(\ref{neg1}), contradicts (\ref{carac_bep}).

Suppose now that $g$ is the solution of $ABEP$. The property that
$|g|=1$ 
on $J$ has been proven in Theorem \ref{exunc}. In order to let $n$
tend to infinity, we rewrite (\ref{carac_fbep}) with self-explaining 
notations as
\begin{equation}
\label{orth}
\langle k_n-f,e^{im\theta} \rangle_{I}+\sum_{j=1}^{r(n)} \lambda^n_j k_n(e^{i\theta_j^n})\overline{ e^{im\theta_j^n}} =0,\quad \forall m \in \{0 \dots n\},\,. 
\end{equation}
We define $(\Lambda_n),\;\;n \in \NN,$ to be a family of linear forms
on $C(J)$ defined as
$$\;\Lambda_n(u)=\sum_{j=1}^{r(n)} \lambda^n_j
k_n(e^{i\theta_j^n})u(e^{\theta_j^n}),\quad \forall u \in C(J).$$
Equation (\ref{lambda_bound}) shows that $(\Lambda_n)$ is a bounded
sequence in the dual $C(J)^*$ which by the Banach-Alaoglu theorem
admits a weak-* converging subsequence whose limit we call
$\Lambda$. Moreover, the Riesz representation theorem ensures the existence of 
a complex measure $\mu$ to represent $\Lambda$ so that, appealing to
Theorem \ref{conv_n} and  taking the limit in (\ref{orth}), we obtain
\begin{equation}
\langle g_0-f,e^{im\theta} \rangle_{I}+\int_{J} \overline{e^{im\theta}}d\mu =0 ,
\quad \forall m \in \NN.
\label{eqFMR}.
\end{equation}
Now, the F. and M. Riesz theorem
asserts that the measure which is $\mu$ on $J$ and $(g_0-f)d\theta$ on
$I$ is absolutely continuous with respect to Lebesgue measure, because its 
Fourier coefficients of nonnegative index do vanish, by \eqref{eqFMR}. 
Therefore there is $v \in L^1(J)$ such that,
$$
\langle g_0-f,e^{im\theta} \rangle_{I}+ \langle v,e^{im\theta} \rangle_{J} =0,
\quad \forall m \in \Bbb N\,,
$$
which is equivalent to
\begin{equation}
\label{first_inf_equ}
\langle g_0-f,e^{im\theta}\rangle _{I}+\langle \lambda g_0,e^{im\theta} \rangle_{J} =0 ,
\quad \forall m \in \Bbb N,
\end{equation}
where we have set $\lambda(z)=v(z)\overline{g_0(z)}\ \forall z \in
J$.
Equation (\ref{first_inf_equ}) means that
$$ {\bf P}_+((g_0-f)\chi_I)=- {\bf P}_+(0 \vee \lambda g_0),$$ 
which indicates that ${\bf P}_+(0 \vee \lambda g_0)$ lies in $H^2$. 
Thus, thanks to Lemma \ref{lemme_integr}, we get that
\begin{equation}
\label{final_orth}
 \langle g_0-f,u \rangle_{I}+\langle \lambda g_0,u \rangle_{J} =0,\quad
\forall u \in H^{2,\infty}_J.
\end{equation}
In order to prove the realness as well as the nonnegativity
of $\lambda$,
we pick $h \in C^\infty_{c,\RR}(I)$, the space of smooth real-valued
functions with compact support on $I$,
and we consider its Riesz-Herglotz transform
\begin{equation}
b(z) = \frac{1}{2 \, \pi} \, \int_{I} \frac{e^{it} + 
z}{e^{it} - z} \, h(e^{it}) \, d t = \frac{1}{2 \, \pi} \,
\int_{\TT} \frac{e^{it} +  
z}{e^{it} - z} \, \chi_I(e^{it}) \, h(e^{it}) \, d t \,. 
\end{equation}
It is standard that $b$ is continuous on
$\overline{\DD}$ \cite[ch. III, thm. 1.3]{Garnett}. 
For $t \in \RR$, define 
$\omega_{t} = \exp (t b)$
which is the outer function whose modulus is equal to 
$\exp t h$ on $I$ and 1 on $J$. 
The function $g_0\, \omega_{\lambda}$ is a candidate approximant
in problem $ABEP$, hence
$t \mapsto \|f - g_0\, \omega_{t}\|_{L^2(I)}^2$ reaches a
minimum at $t= 0$. 
By the boundedness of $b$, we may differentiate this function with respect to
$t$ under the integral sign, and equating the derivative to 0 at
$t=0$ yields
$$0= \mbox{\rm Re} \langle (f-g_0)\overline{g_0},b \rangle_{I}=
\mbox{\rm Re}\langle(f-g_0),bg_0\rangle_{I}.$$ 
In view of \eqref{final_orth}, it implies that
$$0=\mbox{\rm Re}\langle\lambda g_0,bg_0\rangle_{J} =\mbox{\rm Re} \langle \lambda ,b \rangle_{J},
$$
where we used that $|g_0|\equiv1$ on $J$.
Remarking that $b$ is pure imaginary on $J$, this means 
$$  \langle Im(\lambda),b \rangle_{L^2(J)}=0,\quad \forall h \in  C^\infty_{c,\RR}(I).$$ 
Letting $h=h_m$ range over a sequence of smooth positive functions
which are approximate identies,
namely of unit $L^1(I)$-norm and  supported on the arc $[\theta-1/m,\theta+1/m]$
with $e^{i\theta}\in I$, we get in the limit, as $m\to\infty$,  that
$$  \langle Im(\lambda),(e^{i\theta}+.)/(e^{i\theta}-.) \rangle_{J}=0,\quad 
e^{i\theta}\in I.$$ 
Then, appealing to he Plemelj-Sokhotski formulas as
in the proof of Theorem \ref{exunc}, this time on $J$,
we obtain that $Im(\lambda)=0$ which
proves that $\lambda$ is real-valued.
Note that the argument based on the Plemelj-Sokhotski formulas and the 
Hahn-Banach theorem together imply that the space generated by 
$\xi\mapsto (e^{i\theta}+\xi)/(e^{i\theta}-\xi)$, as $e^{i\theta}$ ranges over
an infinite compact subset lying interior to $J$, 
is dense in $L^p(I)$ for $1<p<\infty$.
In fact using the F. and M. Riesz theorem and the Plemelj-Sokhoski formulas, 
it is easy to see that 
such functions are also uniformly dense
in $C(\overline{I})$.
Then, using that $ABEP$ is a convex problem, we obtain 
upon differentiating once more that
$$ \mbox{\rm Re} \langle (g_0-f)\bar{g}_0,b^2 \rangle_{I}\,\, \geq 0, $$ 
which leads us by \eqref{final_orth} to 
$$\mbox{\rm Re} \langle \lambda ,((e^{i\theta}+.)/(e^{i\theta}-.))^2 \rangle_{J}=
 \mbox{\rm Re} \langle \lambda g_0,g_0((e^{i\theta}+.)/(e^{i\theta}-.))^2 \rangle_{J}\, \leq 0,\qquad e^{i\theta}\in I.$$ 
By the density property just mentioned this implies that
$((e^{i\theta}+.)/(e^{i\theta}-.))^2_{|_{\bar I}}$ is dense in the set of 
nonpositive continuous functions on $\overline{I}$, therefore $\lambda\geq0$.
Note also that (\ref{first_inf_equ}) implies $(f-g_0)
\vee \lambda g_0 \in \bar{H}^{1}$, hence it cannot vanish
on a subset of $\TT$ of positive measure unless it is the zero
function. But this would imply $f=g$ a.e on $I$ which contradicts the
hypothesis. This yields $\lambda>0$ a.e on $J$.      
\hfill \boite\\

\section{A numerical example}
\label{numerique}
For practical applications the continuous constraint of $PBEP$ on the arc $J$ is discretized in $m+1$ points. Suppose that $J=\{e^{i t}, t \in [-\theta,\theta]\}$, for some $\theta \in [0,\pi]$. Call $J_m$ the discrete version of the arc $J$ defined by 
$$J_m=\{e^{i t}, t \in \{-\theta+ \frac{2k\theta}{m}, k \in \{0 \dots m\} \} $$ we define following auxiliary extremal problem:
\begin{description}
\item[DBEP(n,m)] \mbox{}\\
For $f \in L^2(I)$, find $k_{n,m} \in P_n$ such that $\forall t \in J_m$
$|k_{n,m}(t)| \leq 1$
and 
\begin{equation} 
\|f - k_{n,m}\|_{L^2(I)} = \min_{\stackrel{g \in P_n}{|g| \leq 1 \,
    \mbox{a.e. on} \, J_m}} \|f - g\|_{L^2(I)} \, .
\label{eq-approxpd}
\end{equation}  
\end{description}

For the discretized problem {\bf DBEP(n,m)}, the following holds.

\begin{theorem}
\label{thmnum}
For $\lambda=(\lambda_0,\dots,\lambda_m) \in \RR^{m+1}$ and $g \in P_n$ define the Lagrangian $$L(\lambda,g)=\|f-g\|_{L^2(I)} + \sum_{k=0}^{m} \lambda_k (|g(e^{i (-\theta+ \frac{2k\theta}{m})}|^2-1)$$, then
\begin{itemize} 
\item Problem {\bf DBEP(n,m)} has a unique solution $k_{n,m}$,
\item $k_{n,m}$ is also the unique solution of the concave maximisation problem:
\begin{equation}
\label{minmax}
\mbox{\rm to find $g_{opt}$ and $\lambda_{opt}$ solving for  }\quad
\max_{\lambda \geq 0} \min_{g \in P_n} L( \lambda,g),
\end{equation}
 where $\lambda \geq 0$ means that each component of $\lambda$ is non negative.
\item For a fixed $n$, $\lim_{m \rightarrow \infty} k_{n,m}=k_n$ in 
$P_n$.
\end{itemize}
\end{theorem}
The proof of Theorem \ref{thmnum} follows from standard convex optimization
theory, using in addition that the $sup$-norm of the derivative of a
polynomial of degree $n$ on $\TT$ is controlled by the values it assumes 
at a set of $n+1$ points. This depends on Bernstein's inequality and
on the argument using Lagrange interpolation polynomials used in the Remark 
after Theorem \ref{exunc}.

In the minmax problem (\ref{minmax}) , the minimization is a quadratic convex problem. It can be tackled efficiently by solving the
critical point equation which is a linear system of equations similar to (\ref{carac_fbep}). Eventually,  an explicit expression of the gradient and 
of the hessian of the concave maximization problem (\ref{minmax}) 
allows us for a fast converging computational procedure 
to estimate $k_{n,m}$. 

\begin{figure}
\begin{center}
\includegraphics[width=0.9\linewidth]{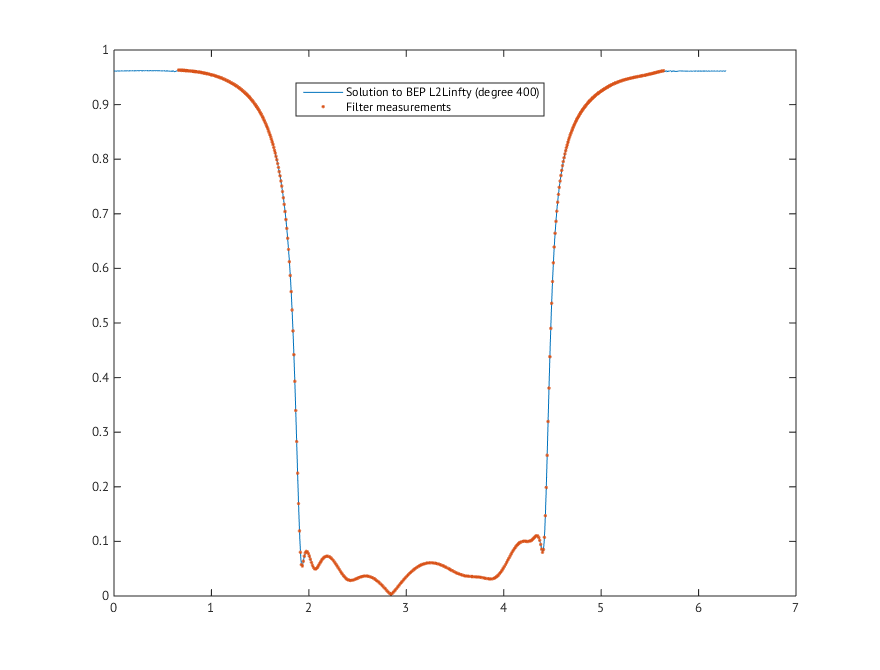}
\end{center}
\caption{Solution of DBEP at hand of partial scattering measurements of a microwave filter}  
\label{cnes}
\end{figure}

Figure (\ref{cnes}) represents a solution to problem {\bf DBEP(n,m)}, where $f$ is obtained from partial measurement of the scattering reflexion parameter of a wave-guide microwave filter by the CNES (French Space Agency). The problem is solved for $n=400$ and $m=800$, while the constraint on $J$ has been renormalized to $0.96$ (instead of 1). 
The modulus of $k_{400,800}$ is plotted as a blue continuous line while the measurements $|f|$ appear as red dots.  As the reader can see, the fit is 
extremely good.

\bibliographystyle{plain}
\bibliography{biblio2001}
\end{document}